\documentclass[12pt,a4paper,titlepage]{article}

\usepackage{amsmath,amssymb,amsthm}
\usepackage{graphicx}

\newtheorem{theorem}{Theorem}
\newtheorem{lemma}[theorem]{Lemma}

\usepackage{amsmath,amssymb,amsfonts,latexsym,epsfig,graphicx,longtable,wrapfig}

\begin{document}

\begin{center}
{\Large
Maximum induced forests in random graphs}\\
\end{center}

\begin{center}
{\large M. Krivoshapko\footnotemark[1], M. Zhukovskii\footnote[1]{Moscow Institute of Physics and Technology (National Research University), Department of Discrete Mathematics, Dolgoprudny, Moscow Region, Russian Federation.}\footnote[2]{Adyghe State University, Caucasus mathematical center, Maykop, Republic of Adygea, Russian Federation; The Russian Presidential Academy of National Economy and Public Administration, Moscow, Russian Federation; Moscow Center for Fundamental and Applied Mathematics, Moscow, Russian Federation.\\ zhukmax@gmail.com}\\}
\end{center}

\vspace{0.5cm}

\begin{center}
Abstract\\
 
We prove that with high probability maximum sizes of induced forests in dense binomial random graphs are concentrated in two consecutive values.\\ 
\end{center}

	\section{Introduction}
	
	Given a graph $G$, its subgraph $H$ is called {\it induced}, if any two vertices $u,v$ of $H$ are adjacent in $H$ if and only if they are adjacent in $G$. {\it The independence number} of $G$ is the maximum number of vertices in an induced subgraph of $G$ that contains no edges. {\it A forest} is an acyclic graph. Everywhere below, {\it the size} of $G$ is the number of vertices in $G$.
 	
	It is very well known (see~\cite{ErdosBol,Grimmett,Matula,Matula2}) that the independence number of the binomial random graph $G(n,p=\mathrm{const})$ (in this graph, every pair of distinct vertices from $\{1,\ldots,n\}$ is adjacent with probability $p$ independently of the others) is concentrated in two consecutive values. In other words, there is a function $f(n)$ such that with high probability (i.e., with probability tending to 1 as $n\to\infty$) the independence number of $G(n,p)$ equals either $f(n)$ or $f(n)+1$. In~\cite{Dutta}, the same concentration result was obtained for the maximum size of an induced path in $G(n,p)$ and for the maximum size of an induced cycle in $G(n,p)$. Finally, in~\cite{Zhuk}, the following 2-point concentration result was obtained for the maximum size of an induced tree.

\begin{theorem}[\cite{Zhuk}]
There exists an $\varepsilon>0$ such that with high probability the maximum size of an induced tree in $G(n,p)$ equals either $\left\lfloor2\log_{1/(1-p)}(enp)+2+\varepsilon\right\rfloor$ or $\left\lfloor2\log_{1/(1-p)}(enp)+3+\varepsilon\right\rfloor$.
\label{th_trees}
\end{theorem}

In this paper, we prove that the same holds true for the maximum size of an induced forest in $G(n,p)$.

\section{The result}

\begin{theorem}
There exists an $\varepsilon>0$ such that with high probability the maximum size of an induced forest in $G(n,p)$ equals either $\left\lfloor2\log_{1/(1-p)}(enp)+2+\varepsilon\right\rfloor$ or $\left\lfloor2\log_{1/(1-p)}(enp)+3+\varepsilon\right\rfloor$.
\label{main}
\end{theorem}

Let $\varepsilon>0$ be the constant from Theorem~\ref{th_trees}. Let $T_n$, $F_n$ be the maximum sizes of an induced tree and an induced forest in $G(n,p)$ respectively. Since a tree is a forest itself, we get that $F_n\geq T_n$. Then, the lower bound in Theorem~\ref{main} follows from Theorem~\ref{th_trees}. The upper bound of Theorem~\ref{th_trees} follows from Markov's inequality. More formally, let $X_n$ be the number of induced trees in $G(n,p)$ of size $\left\lfloor2\log_{1/(1-p)}(enp)+4+\varepsilon\right\rfloor$. In~\cite{Zhuk}, it is proven that ${\sf E}X_n\to 0$. By Markov's inequality, it immediately implies that ${\sf P}(X_n=0)\to 1$. Now, let $Y_n$ be the number of induced forests in $G(n,p)$ of size $\left\lfloor2\log_{1/(1-p)}(enp)+4+\varepsilon\right\rfloor$. In Section~\ref{proof}, we prove the following lemma.

\begin{lemma}
There exists $C>0$ such that ${\sf E}Y_n\leq C{\sf E}X_n$.
\label{expectation_bound}
\end{lemma}

Lemma~\ref{expectation_bound} and the above arguments imply that ${\sf P}(Y_n>0)\leq{\sf E}Y_n\to 0$ as $n\to\infty.$ Theorem~\ref{main} follows.

%Из того, что асимптотика для математического ожидания количества индуцированных лесов на $k$ вершинах такая же, как для математического ожидания количества индуцированных деревьев на $k$ вершинах следует, что:
%$$ \text{E} Y_{\lceil  \hat{k} + \varepsilon \rceil} < C\cdot \text{E} X_{\lceil  \hat{k} + \varepsilon \rceil}  \to 0.$$

%Отсюда и из неравенства Маркова получаем $P(Y_{\lceil  \hat{k} + \varepsilon \rceil} < 1) \to 1$. Таким образом, для размера наибольшего леса также верен результат о концентрации значений в 2 точках.

\section{Proof of Lemma~\ref{expectation_bound}}
\label{proof}

%Покажем, что асимптотика математического ожидания числа индуцированных лесов на $k$ вершинах в случайном графе $G(n, p)$ совпадает с асимптотикой числа деревьев, а именно докажем лемму 1.

%Сформулируем теорему Таннери \cite{Tannery}, которая понадобится нам при доказательстве.

%\begin{theorem}[Теорема Таннери]
%Пусть $S_{n}=\sum\limits_{k=0}^{\infty }a_{k}(n)$ и выполнены условия:
%\begin{enumerate}
%    \item $\lim \limits_{n \rightarrow \infty} a_{k}(n)=b_{k}$,
%    \item для всех $k$ сущетвует $M_k$, т.ч. $\forall n \ \ a_{k}(n) < M_k$,
%    \item для $M_k$ из предыдущего условия $\sum\limits_1^{\infty} M_k \ \ < \infty$.
%\end{enumerate}
%Тогда предел $\lim\limits_{n \rightarrow \infty} S_n$ существует и равен $\sum\limits_1^{\infty} b_k$.
%\end{theorem}

Set $K=\left\lfloor2\log_{1/(1-p)}(enp)+4+\varepsilon\right\rfloor$.  For $\ell\in\{1,\ldots,K\}$, let  $Y_{n, \ell}$ be the number of induced forests in $G(n,p)$ on $K$ vertices with $\ell$ connected components. For $k\geq 2$ and $\ell\in\{1,\ldots,k\}$, let $\varphi_{\ell}(k)$ be the number of forests with $\ell$ components on a labeled set of $k$ vertices. Then
$$
{\sf E}Y_n = \sum_{\ell = 1}^{K} {\sf E} Y_{n,\ell}= \sum_{\ell = 1}^{K} {n\choose K}\varphi_{\ell}(K) p^{K - l}(1 - p)^{{K\choose 2} - K + l}= 
$$
\begin{equation}
{n\choose K}p^{K - 1}(1 - p)^{{K\choose 2} - K + 1}\sum_{\ell = 1}^{K} \varphi_{\ell}(K) p^{-\ell + 1}(1 - p)^{\ell - 1}=
{\sf E}X_n\sum_{\ell = 1}^{K}g_{\ell}(K),
\label{E_compute}
\end{equation}
where
$$
g_{\ell}(k) = \frac{\varphi_{\ell}(k) \cdot ((1 - p)/p)^{\ell - 1}}{k^{k - 2}}.
$$
Set $g_{\ell}(k)=0$ for all $\ell>k$.\\

Since (see~\cite[Section 4.3]{Moon}) $\lim_{k \to +\infty}\frac{\varphi_{\ell}(k)}{k^{k - 2}} = \frac{(1/2) ^ {\ell- 1}}{(\ell -1)!}$, we get that $\lim_{k\to\infty}g_{\ell}(k)=\frac{([(1-p)/(2p)]^ {\ell- 1}}{(\ell -1)!}$. Then, due to Tannery's theorem~\cite[Theorem 3.30]{Tannery}, to prove that 
\begin{equation}
\lim_{k\to\infty}\sum_{\ell = 1}^{k}g_{\ell}(k)=\sum_{\ell = 1}^{\infty}\lim_{k\to\infty}g_{\ell}(k)=\sum_{\ell=1}^{\infty}\frac{([(1-p)/(2p)]^ {\ell- 1}}{(\ell -1)!}<\infty,
\label{finiteness}
\end{equation}
it is sufficient to show that, for every $\ell\in\mathbb{N}$, there exists $M_{\ell}$ such that, for all $k\geq 2$, $g_{\ell}(k) < M_{\ell}$ and $\sum_{\ell=1}^{\infty}M_{\ell}<\infty.$ Notice that, from~(\ref{E_compute})~and~(\ref{finiteness}), Lemma~\ref{expectation_bound} follows.\\

%$$ \lim_{k \to +\infty}\frac{\text{E} Y_{k,l}}{ k^{k-2} \cdot C^k_n \cdot p^{k - 1} \cdot (1 - p)^{C^2_k - k + 1}} = \frac{(0.5  \cdot (1-p)/p) ^ {l- 1} }{(l -1)!}.$$
%В формулировке теоремы Таннери положим $$b_l = \frac{(0.5  \cdot (1-p)/p) ^ {l- 1} }{(l -1)!}.$$

Let $\ell\geq 2$, $m\in\{1,\ldots,k-1\}$. Since the number of forests on $\{1,\ldots,k\}$ with $\ell$ components such that the component containing the vertex $k$ has exactly $k-m$ vertices equals ${k-1\choose m}(k - m)^{k - m - 2}\varphi_{\ell - 1}(m)$, we get the following (recall that $\varphi_{\ell-1}(m)=0$ when $m<\ell-1$):
$$
\varphi_{\ell}(k) = \sum_{m = \ell-1}^{k - 1} {k-1\choose m} (k - m)^{k - m - 2} \varphi_{\ell - 1}(m).
$$
Therefore,
\begin{equation} 
g_{\ell}(k) = \frac{1 - p}{p} \sum\limits_{m=\ell-1}^{k - 1} {k - 1\choose m}\frac{m^{m - 2}(k - m)^{k - m - 2}}{k^{k - 2}} g_{\ell - 1}(m) .
\label{recur}
\end{equation}

Let $M_{\ell} = \max_{k\in\mathbb{N}} g_{\ell}(k)$ (the maximum exists since $g_{\ell}(k)$ has a finite limit as $k\to\infty$). Then, (\ref{recur}) implies that
\begin{equation}
g_{\ell}(k)\leq\frac{1 - p}{p}M_{\ell-1} \sum\limits_{m=\ell-1}^{k - 1} {k - 1\choose m} \frac{m^{m - 2}(k - m)^{k - m - 2}}{k^{k - 2}}=\frac{1 - p}{p}M_{\ell-1} \sum\limits_{m=\ell-1}^{k - 1}f(m,k),
\label{g_upper}
\end{equation}

%Notice that $f(m,k):={k - 1\choose m} m^{m - 2}(k - m)^{k - m - 2}$ is exactly the number of forests consisting of 2 components on $\{1,\ldots,k\}$ such that the component containing the vertex $k$ has $k-m$ vertices. Then, the total number of forests consisting of 2 components equals 
%$$
%\sum_{m=1}^{k - 1}  f(m,k)=\frac{1}{2}k^{k-4}(k-1)(k+6)=\frac{1}{2}k^{k-2}(1+o(1))
%$$  
%(see~\cite{Moon}). Therefore, there exists a finite $\max_{k\in\mathbb{N}}\sum_{m=\ell-1}^{k - 1} \frac{f(m,k)}{k^{k-2}}=:c_{\ell}$. 
where 
$$
f(m,k)={k - 1\choose m}\frac{m^{m - 2}(k - m)^{k - m - 2}}{k^{k-2}}.
$$
Let us show that there exists $C>0$ such that $\sum_{m=\ell-1}^{k - 1} f(m,k)\leq \frac{C}{\ell}$ for all $k$ and $\ell$. If the latter is true, then, by~(\ref{g_upper}), we get that $M_{\ell}\leq\frac{1-p}{p}\frac{C}{\ell}M_{\ell-1}$ that immediately implies the desired finiteness of $\sum_{\ell=1}^{\infty}M_{\ell}$.\\

By Stirling's approximation $\sqrt{2\pi n}(n/e)^n<n!<\sqrt{2\pi n}(n/e)^n e^{\frac{1}{12}n}$ (\cite{stirling}), we get that, for all $k\geq 3$ and $m\leq k-2$,
$$
  f(m,k)<\frac{\sqrt{2\pi(k-1)}(k-1)^{k-1} e^{\frac{1}{12(k-1)}}}{2\pi\sqrt{m(k-1-m)}m^m(k-1-m)^{k-1-m}}m^{m - 2}\frac{(k - m)^{k - m - 2}}{k^{k-2}}=
$$
$$
\sqrt{\frac{(k-1)}{2\pi(k-1-m)}}\frac{k^2}{(k-m)(k-1)}\left(1-\frac{1}{k}\right)^{k}\left(1+\frac{1}{k-m-1}\right)^{k-m-1}\frac{e^{\frac{1}{12(k-1)}}}{m^2\sqrt{m}}<
$$
$$
 c\frac{k^{3/2}}{m^{5/2}(k-m)^{3/2}}
$$
for some constant $c>0$. If $m=k-1$, then $f(m,k)=\frac{(k-1)^{k-3}}{k^{k-2}}<\frac{1}{k}$. So, the above bound is also true in this case.

Notice that the function $v(x) =  \frac{k^{3/2}}{x^{5/2}(k - x) ^ {3/2}}$ is convex on $(0,k)$. Therefore, for all $k\geq 2$ and $2\leq\ell\leq k$,
$$ 
\sum_{m = \ell - 1}^{k - 1} \frac{k^{3/2}}{m^{5/2}(k -m) ^ {3/2}} \leq v(k - 1) + v(\ell - 1) + \int_{\ell - 1}^{k - 1} v(x) dx=
$$
$$
\frac{k^{3/2}}{(k-1)^{5/2}}+\frac{k^{3/2}}{(\ell-1)^{5/2}(k - \ell+1) ^ {3/2}}-\left.\frac{2 (k^2 + 4 k x - 8 x^2)}{3 k^{3/2}x^{3/2} \sqrt{k - x} }\right|_{\ell-1}^{k-1}<
$$
$$
\frac{6}{k}+\frac{3}{\ell-1}+\frac{2 (3k^2-12k+8)}{3 k^{3/2}(k-1)^{3/2}}+\frac{2 (k^2 + 4 k(\ell-1) - 8(\ell-1)^2)}{3 k^{3/2}(\ell-1)^{3/2}\sqrt{k - \ell+1}}.
$$
Since 
$$
\frac{2 (3k^2-12k+8)}{3 k^{3/2}(k-1)^{3/2}}<\frac{6(k-2)^2}{3 k^{3/2}(k-1)^{3/2}}<\frac{2\sqrt{k-1}}{3k\sqrt{k}}<\frac{2}{3k}
$$
and
$$
 \frac{2 (k^2 + 4 k(\ell-1) - 8(\ell-1)^2)}{3 k^{3/2}(\ell-1)^{3/2}\sqrt{k - \ell+1}}<
 \frac{2[(k+2(\ell-1))^2-16(\ell-1)^2]}{3 k^{3/2}(\ell-1)\sqrt{(k -(\ell-1))(\ell-1)}}<
$$
$$
 \frac{2(k-2(\ell-1))(k+6(\ell-1))}{k^2(\ell-1)}<\frac{2(k+6k)}{k(\ell-1)}=\frac{14}{\ell-1},
$$
the desired bound follows.

\section{Discussions}

We have proved that with high probability the maximum size of an induced forest in $G(n,p)$ is concentrated in 2 consecutive points: $\left\lfloor2\log_{1/(1-p)}(enp)+2+\varepsilon\right\rfloor$ and $\left\lfloor2\log_{1/(1-p)}(enp)+3+\varepsilon\right\rfloor$. Moreover, from our arguments, it follows that the difference between it and the maximum size of an induced tree is not bigger than 1.
It would be interesting to extend these results for $p=o(1)$. In sparse settings (particularly, for $p=c/n$ and $p=c\ln n/n$), some bounds are known~\cite{Vega,Luczak,Palka} for the maximum size of an induced tree, but even the exact asymptotics of these maximum sizes are unknown.

\section{Acknowledgements}

The reported study was funded by RFBR according to the research project N 20-04-60524.

\end{document}